\newif\ifNDofficial
\NDofficialfalse

\ifNDofficial
\documentclass[11pt,%
final,%
tightenlines,%
notitlepage,%
twoside,%
onecolumn,%
nobibnotes,%
nofootinbib,%
superscriptaddress
]{rcd-nd-eng}
\usepackage{epsfig}
\usepackage{graphicx}
\usepackage{wrapfig}
\usepackage{psfrag}
\usepackage{substr}
\usepackage{verbatim}
\usepackage[tight]{minitoc}
\usepackage{subfigure}
\usepackage{longtable}
\usepackage[all]{xy}
\usepackage{xcite}
\usepackage{xr}
\usepackage[colorlinks,linkcolor=blue,pagecolor=blue,citecolor=blue,urlcolor=blue,hypertexnames=false]{hyperref}
\else
\documentclass[11pt,twoside]{article}
\usepackage[T1]{fontenc}
\usepackage[utf8]{inputenc}
\usepackage[english]{babel}
\usepackage[a4paper,left=2.55cm,right=2.55cm,top=2.15cm,bottom=2.35cm]{geometry}
\usepackage[colorlinks=true,linkcolor=blue!55!black,citecolor=blue!55!black,urlcolor=blue!55!black,hypertexnames=false]{hyperref}
\usepackage{fancyhdr}
\usepackage{xcolor}
\usepackage{microtype}
\newcommand{\keywords}[1]{\par\medskip\noindent\textbf{Keywords:} #1\par}
\fi

\usepackage{amsmath,amssymb}
\usepackage{graphicx}
\usepackage{booktabs}
\usepackage{array}
\usepackage{enumitem}
\usepackage{algorithm}
\usepackage{algpseudocode}
\graphicspath{{figs/}}
\allowdisplaybreaks

\def\R{\mathbb{R}}
\def\Expct{\mathbb{E}}
\def\Zset{\mathcal{Z}}
\def\Gap{\mathop{\rm Gap}\nolimits}
\def\supp{\mathop{\rm supp}\nolimits}
\def\inner#1#2{\langle #1,#2\rangle}
\def\norm#1{\left\|#1\right\|}
\def\hG{\widehat{G}}
\def\hz{\hat z}
\newenvironment{pf}[1][Proof]{\par\noindent\textit{#1.}\ }{\hfill$\square$\par\medskip}

\ifNDofficial

\rcdthmsep=\medskipamount
\newrcdthm{teo}{Theorem}
\newrcdthm{lem}{Lemma}
\newrcdthm{pro}{Proposition}
\newrcdthm{assum}{Assumption}
\newrcdthm{cor}{Corollary}
\newrcdthm{con}{Conjecture}
\newrcdthm{fed}{Definition}
\newrcdthm[remark]{rem}{Remark}
\newrcdthm[example]{exl}{Example}
\newrcdthm[example]{com}{Comment}
\else
\usepackage{amsthm}
\newtheorem{teo}{Theorem}[section]
\newtheorem{lem}[teo]{Lemma}
\newtheorem{pro}[teo]{Proposition}
\newtheorem{assum}[teo]{Assumption}
\newtheorem{cor}[teo]{Corollary}

\theoremstyle{definition}
\newtheorem{fed}[teo]{Definition}

\theoremstyle{remark}

\fi

\hypersetup{%
  pdftitle={Federated Affine Variational Inequalities: Exact Averaging and Complexity Limits},%
  pdfauthor={Darina M. Dvinskikh, Kirill A. Novoselov, and Alexander V. Gasnikov},%
  pdfsubject={Exact averaging and communication-locality limits for federated affine variational inequalities},%
  pdfkeywords={federated learning, intermittent communication, stochastic variational inequalities, local stochastic extragradient, affine monotone operators, exact averaging, communication complexity}%
}

\begin{document}

\ifNDofficial
\def\mtctitle{Contents}%
\def\mtcSfont{\small\rmfamily\upshape\mdseries}
\makeatletter
\let\mtc@rule=\relax
\makeatother
\allowdisplaybreaks

\widowpenalty =10000
\clubpenalty = 10000

\setcounter{minitocdepth}{2}

\renewcommand{\hm}[1]{#1\nobreak\discretionary{}%
{\hbox{$\mathsurround=0pt #1$}}{}}

\IfFileExists{letterspacing.tex}{%
  \let\NDAlgorithmicComment\Comment
  \let\Comment\Undefined
  \let\EndComment\Undefined
  \input{letterspacing.tex}%
  \let\Comment\NDAlgorithmicComment
  \let\EndComment\relax
  \def\trk#1#2{\mbox{}\letterspace to #1\naturalwidth{#2}}%
}{%
  \providecommand{\trk}[2]{\mbox{#2}}%
}

\renewcommand\refname{References}
\renewcommand{\figurename}{Fig.}

\JournalName{Russian Journal of Nonlinear Dynamics}
\JournalNameEng{Rus. J. Nonlin. Dyn.}
\OrigYearOfIssue{2026}
\OrigVolumeNo{}
\OrigIssueNo{}
\NDcite{https://nd.ics.org.ru/}

\setcounter{page}{1}
\Rubrika{\relax}
\CRubrika{\relax}
\SubRubrika{\relax}

\beginpaper
\amsmsc{90C25, 90C15, 65K15, 68Q17, 90C47}
\titlerunning{Exact averaging and complexity limits}
\title{Federated Affine Variational Inequalities: Exact Averaging and Complexity Limits}
\toctitleeng{Federated Affine Variational Inequalities: Exact Averaging and Complexity Limits}
\authorrunning{D.\,M.\,Dvinskikh, K.\,A.\,Novoselov, A.\,V.\,Gasnikov}
\author{D.\,M.\,Dvinskikh}
\tocauthoreng{D.\,M.\,Dvinskikh, K.\,A.\,Novoselov, A.\,V.\,Gasnikov}
\info{Darina M. Dvinskikh}{dmdvinskikh@hse.ru}
{HSE University, Pokrovsky Blvd 11, Moscow, 109028 Russia}
\author{K.\,A.\,Novoselov}
\info{Kirill A. Novoselov}{k.novoselov@innopolis.university}
{Innopolis University, Universitetskaya 1, Innopolis, 420500 Russia}
\author{A.\,V.\,Gasnikov}
\info{Alexander V. Gasnikov}{gasnikov@yandex.ru}
{Innopolis University, Universitetskaya 1, Innopolis, 420500 Russia;\\Steklov Mathematical Institute of the Russian Academy of Sciences, Ulitsa Gubkina 8, Moscow, 119333 Russia;\\Ivannikov Institute for System Programming of the Russian Academy of Sciences, Ulitsa Aleksandra Solzhenitsyna 25, Moscow, 109004 Russia}
\makeatletter
\gdef\@received{}
\gdef\@receivedeng{}
\gdef\@accepted{}
\gdef\@acceptedeng{}
\makeatother
\grants{The work was supported by the Ministry of Science and Higher Education of the Russian Federation, Grant No.~075-15-2024-529.}
\else
\thispagestyle{empty}

\begin{center}
\vspace{0.9cm}
{ \LARGE\bfseries Federated Affine Variational Inequalities:\par
Exact Averaging and Complexity Limits\par}
\vspace{0.45cm}
{\large\bfseries  D.\,M. Dvinskikh, K.\,A. Novoselov, A.\,V. Gasnikov\par}
\end{center}
\fi

\begin{abstract}
We study the classical Local Stochastic Extragradient method (Local SEG) for
monotone variational inequalities under intermittent communication, where
devices perform several local oracle calls between synchronization rounds. For
affine monotone mean operators, including nonsymmetric bilinear saddle systems,
we prove that equal device averaging commutes pathwise with both stages of the
stochastic extragradient update. Local trajectories may therefore drift apart
while their mean follows a centralized stochastic-extragradient recursion
driven by averaged noise. The identity remains valid when individual oracle
realizations are nonlinear and query dependent.

Under conditionally centered oracle errors that are conditionally independent
across devices within each stage, the reduction yields a restricted weak-gap
rate that uses the local sequential depth for deterministic progress and the
aggregate oracle budget across devices for variance reduction. When every
device makes at least two local oracle calls per communication round, the bound
contains no additional mixed client-drift term and is also attained by a
terminal-only aggregation schedule. An affine Gaussian lower bound matches
the depth and sample dependences before statistical saturation for distributed
zero-respecting algorithms; its statistical component applies to arbitrary
randomized algorithms. A complementary lower bound for feasible distributed
zero-respecting methods shows that bounded geometry and query-dependent noise
can create an additional communication-locality obstruction. These results
delineate the exact-averaging regime of unconstrained affine Local SEG and
a complexity boundary that arises under active feasibility constraints.
\end{abstract}


\ifNDofficial
\maketitle\label{article_begin}
{
\else
\par\medskip
\begingroup\raggedright
\noindent\textit{D. M. Dvinskikh:}\\
HSE University, Pokrovsky Blvd 11, Moscow, 109028 Russia;\\
e-mail: \texttt{dmdvinskikh@hse.ru}.\\[1mm]
\textit{K. A. Novoselov:}\\
Innopolis University, Universitetskaya 1, Innopolis, 420500 Russia;\\
e-mail: \texttt{k.novoselov@innopolis.university}.\\[1mm]
\textit{A. V. Gasnikov:}\\
Innopolis University, Universitetskaya 1, Innopolis, 420500 Russia;
Steklov Mathematical Institute of the Russian Academy of Sciences,
Ulitsa Gubkina 8, Moscow, 119333 Russia;\\
e-mail: \texttt{gasnikov@yandex.ru}.

\par\smallskip \small
\noindent\textit{Funding.} This work was supported by Ministry of Science and Higher Education of
the Russian Federation (Grant No. 075-15-2024-529)
\endgroup
\par\bigskip
\hrule\bigskip
\setcounter{page}{1}
\pagestyle{fancy}
\fi

\section{Introduction}
\label{sec:intro}

Variational inequalities (VIs) are a standard framework for constrained convex
optimization, convex--concave saddle problems, and equilibrium models. A VI
asks for a feasible point whose operator value forms a nonnegative inner
product with every feasible displacement. The operator may combine gradient
and rotational components, which extends the model beyond ordinary
minimization.
Extragradient and Mirror-Prox methods address this difficulty by evaluating the
operator twice per iteration
\cite{Korpelevich1976,Nemirovski2004,JNT2011}.

Large-scale stochastic VIs are commonly addressed on distributed systems. Several
devices can collect oracle samples in parallel, whereas synchronization is
often slower than local computation. This motivates the
intermittent-communication (IC) model: each communication opportunity is separated
by several sequential oracle calls on every device
\cite{ArjevaniShamir2015,WoodworthGraph2018,WBSS2021}. The model is a standard
abstraction for clusters and homogeneous federated systems. It distinguishes
two computational resources: the length of the adaptive local trajectory and
the total number of samples collected across all devices.

A central analytical question is whether a method can exploit both resources
simultaneously. Minibatch methods pool parallel samples, so their deterministic
horizon follows the synchronization count. A single-device method traverses
the full local trajectory with one device's sample budget. Local SGD and FedAvg seek to
combine the two advantages in stochastic minimization
\cite{McMahan2017,Stich2019local,Haddadpour2019,Kairouz2021,PatelEtAl2024}.
The analogous question for VIs is more delicate because the second
extragradient query depends on the first and disagreement can enter both stages
of every update.

General intermittent-communication problems therefore exhibit a tension
between parallel sampling and sequential local progress
\cite{WoodworthGraph2018,WBSS2021}. Homogeneous quadratic minimization is an
important exception: linear updates commute with device averaging, and one
terminal communication can attain the optimal rate
\cite{Woodworth2020local,WBSS2021}. We ask whether this exception extends
beyond gradient fields to nonsymmetric affine monotone VIs and the two-stage
extragradient recursion.

We analyze the established Local Stochastic Extragradient method (Local SEG)
in this setting \cite{BeznosikovLocal2021}. A particularly close
comparison is the recent work of Wang and Kale on Local Extra SGD with an
affine mean operator
\cite{WangKale2026}. Their guarantee uses local computation and parallel
sampling, but contains an additional mixed term arising from the control of
client drift. This raises a structural question: does the device mean inherit
that mixed term, or does affinity yield a simpler average trajectory?

For affine mean operators, the device average follows an exact centralized SEG
recursion even when individual oracle realizations depend nonlinearly on the
query. The resulting two-term weak-gap estimate uses both sequential depth and
parallel sampling, and it can be attained with terminal-only aggregation. The
lower bounds determine its order-optimal scope and identify a separate
communication-locality obstruction on bounded domains.

\paragraph{Contributions.}
The paper makes three analytical contributions to the analysis of Local SEG.
\begin{enumerate}
\item We extend the exact-averaging principle for linear-update methods on
homogeneous quadratic objectives \cite{Woodworth2020local} to nonsymmetric
affine monotone VIs and the two-stage Local SEG recursion. The reduction allows
nonlinear, query-dependent oracle realizations and uses only their adapted
moment bounds. In particular, it covers rotational bilinear saddle operators
with a non-gradient mean field.
\item The reduction removes the mixed client-drift term from the Wang--Kale
affine Local Extra SGD analysis \cite{WangKale2026}. It yields a two-resource
weak-gap guarantee and shows that one terminal aggregation suffices for the
prescribed horizon.
\item Two lower bounds describe both optimality and scope. An affine Gaussian
construction matches the depth and sample dependences in the unconstrained
problem. A
separate bounded affine construction exhibits a communication-locality term
for distributed zero-respecting methods when feasibility constraints and
query-dependent noise are admitted.
\end{enumerate}

\paragraph{Related work and scope.}
Distributed VI methods cover heterogeneous operators, decentralized mixing,
similarity assumptions, variance reduction, and local updates
\cite{Beznosikov2021distributed,Rogozin2021,Beznosikov2022decentralized,
Kovalev2022,Beznosikov2023smooth,BeznosikovPillars2023}.
Polyak--Ruppert averaging studies temporal averages of
stochastic-approximation iterates; our identity concerns averages across
simultaneous device trajectories
\cite{Ruppert1988,PolyakJuditsky1992,Mou2020}. Accelerated federated
minimization methods such as FedAc exploit gradient structure and measure
objective suboptimality \cite{YuanMa2020}. Our analysis uses the weak VI gap
and affine monotone fields, so the two directions address complementary
problem classes.
Sadchikov et al. study how the quadratic exact-averaging phenomenon degrades
beyond quadratic structure for Local SGD \cite{SadchikovEtAl2024}. Their
setting concerns gradient fields and one-stage stochastic-gradient updates.
Our setting adds nonsymmetric affine operators, two coupled SEG stages, a
restricted weak-gap criterion, and affine VI lower bounds.
Interpolation and strong-growth assumptions provide a different route to
reduced stochastic error \cite{Vaswani2019,WoodworthSrebro2021}. Under such
noise models, the effective variance of Local SEG may also depend on
disagreement among local query points, so this regime requires a separate
analysis.
The depth lower bound follows classical oracle-complexity, bilinear-chain, and
distributed graph constructions
\cite{NemirovskyYudin1983,WoodworthGraph2018,OuyangXu2021}; its statistical
component uses a standard two-point testing argument \cite{Tsybakov2009}.
For the communication-locality boundary, we adapt the robust zero-chain idea
of Woodworth et al.~\cite{WBSS2021} to a bounded bilinear VI and compute the
weak gap directly. Our contribution is the affine VI adaptation and the
resulting boundary for exact averaging in the homogeneous model.

The main results concern an unconstrained affine mean operator. Because a
global weak gap on an unbounded domain can be infinite away from a solution,
we evaluate accuracy on a fixed ball around a root. This ball localizes only
the merit function. The algorithm, its iterates, and its output remain
unconstrained. The next section gives the formal definitions and assumptions.

\paragraph{Organization.}
Section~\ref{sec:setting} defines the oracle model and Local SEG.
Section~\ref{sec:vi} proves the exact reduction, its rate, and the matching
affine lower bound. Section~\ref{sec:limits} identifies a constrained
communication-locality boundary.
Section~\ref{sec:numerics} gives a numerical illustration. The appendices
contain the lower-bound constructions, technical SEG estimate, and experimental
details.

\section{Problem setting and Local SEG}
\label{sec:setting}

This section introduces the communication model, the stochastic oracle, the
weak-gap criterion, and the Local SEG recursion used in the main results.

\paragraph{Notation.}
Throughout, $\norm{\cdot}$ is the Euclidean norm for vectors and the induced
operator norm for matrices, and $\inner{\cdot}{\cdot}$ is the Euclidean inner
product.

\subsection{Intermittent communication}

\begin{figure}[h]
    \centering
    \includegraphics[width=0.8\textwidth]{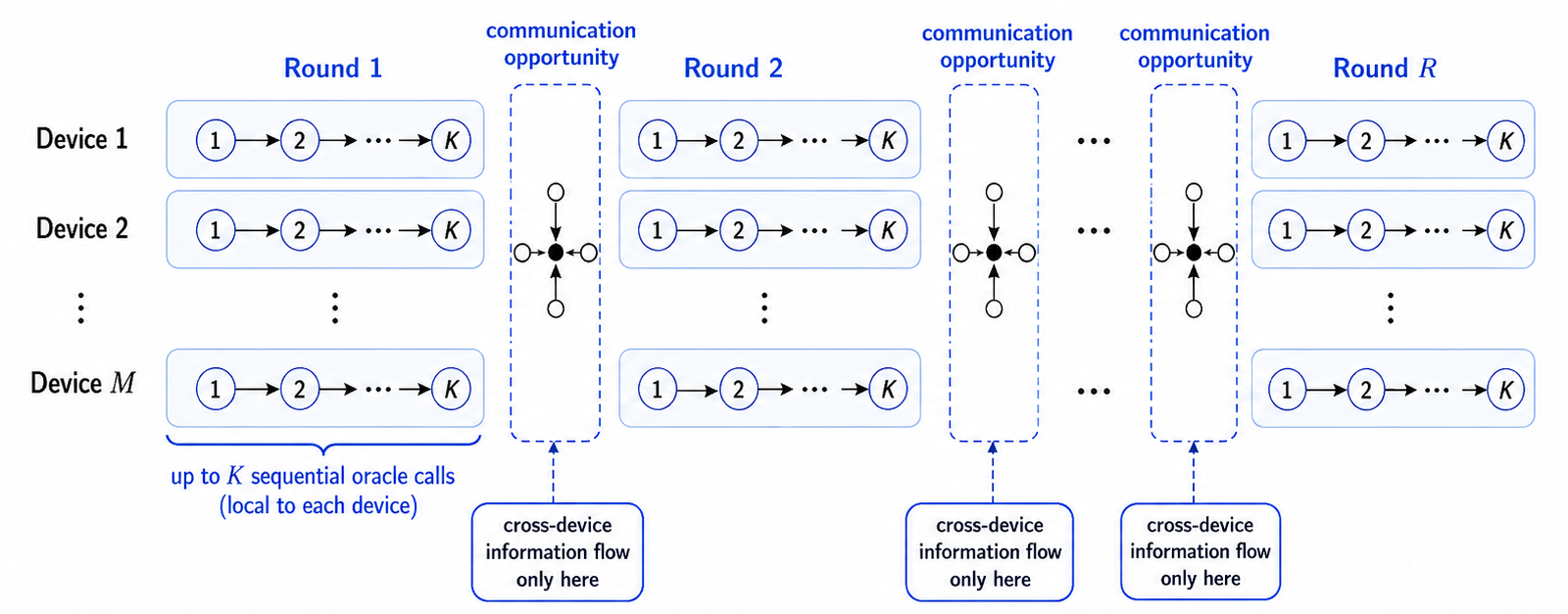}
    \caption{Intermittent-communication model. Each device performs up to $K$
    sequential oracle calls per round, while information can flow across devices
    only at the communication boundaries between rounds.}
    \label{fig:ic-model}
\end{figure}

We use the standard homogeneous intermittent-communication model
\cite{WoodworthGraph2018,WBSS2021}. There are $M$ devices and $R$
local-computation rounds. In each round, every device makes at most $K$
sequential calls to the same stochastic VI oracle, after which communication
between the devices is allowed. Within a round, a device may adapt its queries
to its own preceding oracle replies; information obtained by other devices
becomes available after communication. Thus the total oracle budget is at most
$MKR$, while the maximal sequential depth is $KR$. Figure~\ref{fig:ic-model} illustrates this mechanism.

\paragraph{Lower-bound algorithm class.}
The depth and communication-locality lower bounds concern distributed
zero-respecting IC algorithms. Relative to the standard coordinate basis,
such an algorithm may activate a coordinate only after that coordinate has
appeared in an oracle reply available to it. The formal transcript and support
conditions are given in Appendix~\ref{app:ic-formal}. The statistical lower
bound applies to arbitrary randomized IC algorithms.

\subsection{Variational inequalities, oracle, and merit function}

Let $\Zset\subseteq\R^n$ be nonempty, closed, and convex and let
$G:\Zset\to\R^n$. The variational inequality ($\operatorname{VI}$)
is to find $z^\star\in\Zset$ such that
\begin{equation}
\label{eq:vi}
\inner{G(z^\star)}{z-z^\star}\ge0
\qquad\text{for every }z\in\Zset.
\end{equation}
For continuous monotone $G$, this is equivalent to the weak formulation
$\inner{G(z)}{z^\star-z}\le0$ for every $z\in\Zset$.

We model a deterministic mean operator $G$ accessed through a stochastic
oracle. We therefore state the operator properties and the stochastic-oracle
conditions separately.

When stochastic access is used, let $\Xi$ be a sample space and let
$\hG:\Zset\times\Xi\to\R^n$ denote the oracle response map.

\begin{assum}[Monotonicity]
\label{as:monotone}
The operator $G$ is monotone on $\Zset$: for every $z,z'\in\Zset$,
\[
\inner{G(z)-G(z')}{z-z'}\ge0.
\]
\end{assum}

\begin{assum}[Lipschitz continuity]
\label{as:lipschitz}
For some $L>0$, the operator $G$ is $L$-Lipschitz on $\Zset$: for every
$z,z'\in\Zset$,
\[
\norm{G(z)-G(z')}\le L\norm{z-z'}.
\]
\end{assum}

\begin{assum}[Stochastic oracle]
\label{as:oracle}
For some $\sigma\ge0$, before each call $\mathcal F$ denotes the complete
history and the query $z$ is $\mathcal F$-measurable. A new draw
$\zeta\in\Xi$ satisfies
\[
\Expct\!\left[\hG(z,\zeta)\mid\mathcal F\right]=G(z),\qquad
\Expct\!\left[\norm{\hG(z,\zeta)-G(z)}^2\mid\mathcal F\right]\le\sigma^2.
\]
Randomness used by simultaneous calls is conditionally independent given the
pre-stage history. At a later stage, $\mathcal F$ includes all earlier replies.
\end{assum}

For the unconstrained affine results in Section~\ref{sec:vi},
$\Zset=\R^n$. We call $z^\star$ a \emph{root} of $G$ if
$G(z^\star)=0$; such a root solves~\eqref{eq:vi}. Whenever an affine instance
has a specified root $z^\star$ and radius $\Omega>0$, we measure any candidate
$x\in\R^n$ by the root-centered restricted weak gap
\begin{equation}
\label{eq:gap-omega}
\Gap_\Omega(x)
:=\sup_{\norm{z-z^\star}\le\Omega}\inner{G(z)}{x-z}.
\end{equation}
The comparison ball in~\eqref{eq:gap-omega} only localizes the merit function;
the algorithm, its iterates, and its output remain unconstrained in $\R^n$. This
localization is necessary because the unrestricted weak gap on $\R^n$ may be
infinite away from a solution. Since the comparison ball contains $z^\star$
and $G(z^\star)=0$, we have $\Gap_\Omega(x)\ge0$ for every $x$.
Monotonicity further gives $\Gap_\Omega(z^\star)=0$.

\subsection{Local SEG}
\label{sec:free-local-seg}

We use the Local SEG recursion
\cite{BeznosikovLocal2021,WangKale2026}. Algorithm~\ref{alg:local-seg}
records its two oracle stages and makes the equal-averaging schedule explicit.
The initial states are $\mathcal F_0$-measurable, and all devices use a common
deterministic stepsize.

\begin{algorithm}[H]
\caption{Local SEG with an admissible equal-averaging schedule}
\label{alg:local-seg}
\begin{algorithmic}[1]
\Require States $z_0^1,\ldots,z_0^M$, horizon $T$, stepsize $\gamma>0$,
and a synchronization schedule
\State $s_0^m\gets0$ for every $m=1,\ldots,M$
\For{$t=0,\ldots,T-1$}
  \ForAll{devices $m$ in parallel}
    \State $y_t^{m,-}\gets z_t^m-
      \gamma\widehat G(z_t^m,\zeta_t^m)$
  \EndFor
  \State Set $y_t^m\gets M^{-1}\sum_j y_t^{j,-}$ for every $m$ if a
  half-stage synchronization is scheduled; otherwise set
  $y_t^m\gets y_t^{m,-}$
  \ForAll{devices $m$ in parallel}
    \State $z_{t+1}^{m,-}\gets z_t^m-
      \gamma\widehat G(y_t^m,\zeta_{t+1/2}^m)$
    \State $s_{t+1}^m\gets s_t^m+y_t^m$
  \EndFor
  \State Set $z_{t+1}^m\gets M^{-1}\sum_j z_{t+1}^{j,-}$ for every $m$ if
  a terminal-stage synchronization is scheduled; otherwise set
  $z_{t+1}^m\gets z_{t+1}^{m,-}$
\EndFor
\State \Return $\displaystyle
  \widetilde z_T=\frac1{MT}\sum_{m=1}^M s_T^m$
\end{algorithmic}
\end{algorithm}

An admissible synchronization either leaves the corresponding stage outputs
unchanged or replaces all of them by their equal-weight mean. No projection is
applied. We write
\[
\bar z_t:=\frac1M\sum_{m=1}^M z_t^m,
\qquad
\bar y_t:=\frac1M\sum_{m=1}^M y_t^m.
\]
The two local extragradient stages are standard; the optional stage-wise
averaging makes the schedule-independent mean identity explicit.
Algorithm~\ref{alg:local-seg} records the method, and
Section~\ref{sec:vi} develops its exact-reduction analysis. Classical block
Local SEG performs $\tau=\lfloor K/2\rfloor$ completed SEG iterations per
communication round, so $T=R\tau$. The terminal-only implementation leaves
intermediate communication opportunities unused and aggregates the accumulated
vectors $s_T^m$ at the end.

\subsection{Affine mean subclass}

\begin{assum}[Affine mean structure]
\label{as:affine}
Let $\Zset=\R^n$. The mean operator is affine:
\[
G(z)=Jz+q
\]
for some $J\in\R^{n\times n}$ and $q\in\R^n$. The affine condition applies to
the conditional mean operator; individual oracle realizations may be nonlinear
and query dependent.
\end{assum}

For an affine operator on $\R^n$, Assumption~\ref{as:monotone} is equivalent
to $J+J^\top\succeq0$, whereas Assumption~\ref{as:lipschitz} with constant
$L$ is equivalent to $\norm{J}_{\rm op}\le L$. The class includes quadratic
minimization and bilinear saddle problems. A nonsymmetric matrix $J$ allows a
non-gradient mean field. The pathwise identities below
use the affine representation alone. Monotonicity and Lipschitz continuity
enter the weak-gap estimate through Proposition~\ref{prop:freeeg}, while
Assumption~\ref{as:oracle} supplies the conditional moment bounds.

\section{Main results for unconstrained affine VIs}
\label{sec:vi}

This section establishes the exact device-mean recursion, derives its
restricted weak-gap consequence, and proves a corresponding affine lower
bound.

\subsection{Exact local-to-central equivalence}
\label{sec:affine}

The theorem is structural: affinity makes equal averaging commute with both
stages of the standard local SEG update. Its pathwise part is deterministic;
the subsequent moment bounds invoke the stochastic-oracle assumption.

For each iteration, define the oracle residuals
\[
\varepsilon_t^m
:=\widehat G(z_t^m,\zeta_t^m)-G(z_t^m),
\qquad
\varepsilon_{t+1/2}^m
:=\widehat G(y_t^m,\zeta_{t+1/2}^m)-G(y_t^m),
\]
and let
\[
\bar\varepsilon_s:=\frac1M\sum_{m=1}^M\varepsilon_s^m,
\qquad s\in\{t,t+1/2\}.
\]

\begin{teo}[Exact device-mean recursion]
\label{thm:exact}
Let Assumption~\ref{as:affine} hold and run
Algorithm~\ref{alg:local-seg} under any admissible synchronization schedule.
Then, pathwise for every $t=0,\ldots,T-1$,
\begin{equation}
\label{eq:exact-recursion}
\bar y_t
=\bar z_t-\gamma\bigl(G(\bar z_t)+\bar\varepsilon_t\bigr),
\qquad
\bar z_{t+1}
=\bar z_t-\gamma\bigl(G(\bar y_t)+\bar\varepsilon_{t+1/2}\bigr).
\end{equation}
If Assumption~\ref{as:oracle} also holds and $\mathcal F_s$ denotes the
history before oracle stage $s$, then
\[
\Expct\!\left[\bar\varepsilon_s\mid\mathcal F_s\right]=0,
\qquad
\Expct\!\left[\norm{\bar\varepsilon_s}^{2}\mid\mathcal F_s\right]
\le\frac{\sigma^2}{M},
\qquad s\in\{t,t+1/2\}.
\]
\end{teo}

\begin{pf}
At the first stage, affinity gives
\[
\begin{aligned}
\bar y_t
&=\bar z_t-\frac{\gamma}{M}\sum_{m=1}^M
  \bigl(G(z_t^m)+\varepsilon_t^m\bigr)
  =\bar z_t-\gamma\bigl(G(\bar z_t)+\bar\varepsilon_t\bigr).
\end{aligned}
\]
The same calculation at the second stage gives the second identity in
\eqref{eq:exact-recursion}. Optional equal-weight synchronizations leave both
means unchanged. Under Assumption~\ref{as:oracle}, conditional centering gives
the first-moment identities, and conditional independence eliminates cross
terms; for example,
\[
\Expct\!\left[\norm{\bar\varepsilon_t}^2\mid\mathcal F_t\right]
=\frac1{M^2}\sum_{m=1}^M
\Expct\!\left[\norm{\varepsilon_t^m}^2\mid\mathcal F_t\right]
\le\frac{\sigma^2}{M}.
\]
Conditioning on $\mathcal F_{t+1/2}$ proves the second-stage estimate.
\end{pf}

Thus the device mean follows a centralized stochastic-extragradient
recursion with conditionally centered noise whose variance is reduced by the
number of devices. The pathwise identity uses affinity alone; the stochastic
oracle assumption supplies the moment bounds.

\paragraph{Relation to quadratic minimization.}
When $J=J^\top\succeq0$, the operator is the gradient of a quadratic
objective, and the theorem specializes to the exact-averaging principle for
linear-update methods on homogeneous quadratics
\cite{Woodworth2020local}. For nonsymmetric $J$, the same identity applies to
non-gradient operators and preserves the two-stage mean recursion.

\subsection{Weak-gap rate and terminal-only aggregation}

The exact recursion transfers an unprojected SEG estimate with variance
divided by $M$. The common stepsize controls the balance between deterministic
progress and stochastic error.

\begin{cor}[Restricted weak-gap rate]
\label{cor:localseg}
Suppose Assumptions~\ref{as:monotone}, \ref{as:lipschitz},
\ref{as:oracle}, and~\ref{as:affine} hold, and let $z^\star$ be a root of $G$.
Suppose that the initial device states are $\mathcal F_0$-measurable and
$\norm{\bar z_0-z^\star}\le\Omega$
almost surely for some $\Omega>0$. Then, for every admissible synchronization
schedule, every integer $T\ge1$, and every $0<\gamma\le1/(2L)$,
\begin{equation}
\label{eq:localseg-general}
\Expct\!\left[\Gap_\Omega(\widetilde z_T)\right]
\le
\frac{3\Omega^2}{\gamma T}
+\frac{9\gamma\sigma^2}{M}.
\end{equation}
In particular, choose
\[
\gamma
=\min\left\{
\frac{1}{2L},
\frac{\Omega\sqrt M}{\sigma\sqrt T}
\right\},
\]
where the second quantity is interpreted as $+\infty$ when $\sigma=0$.
Then
\begin{equation}
\label{eq:localseg}
\Expct\!\left[\Gap_\Omega(\widetilde z_T)\right]
\le
12\left(
\frac{L\Omega^2}{T}
+\frac{\sigma\Omega}{\sqrt{MT}}
\right).
\end{equation}
\end{cor}

\begin{pf}
Assumptions~\ref{as:monotone} and~\ref{as:lipschitz} provide the regularity
required by Proposition~\ref{prop:freeeg}. Theorem~\ref{thm:exact} gives its
adapted-error recursion with $v=\sigma^2/M$, proving
\eqref{eq:localseg-general}. The restriction $\gamma\le1/(2L)$ ensures
stability. For the stated stepsize,
\[
\frac1\gamma
\le
2L+\frac{\sigma\sqrt T}{\Omega\sqrt M},
\qquad
\gamma\le\frac{\Omega\sqrt M}{\sigma\sqrt T},
\]
with the usual convention when $\sigma=0$. Hence
\[
\frac{3\Omega^2}{\gamma T}
\le
\frac{6L\Omega^2}{T}
+\frac{3\sigma\Omega}{\sqrt{MT}},
\qquad
\frac{9\gamma\sigma^2}{M}
\le
\frac{9\sigma\Omega}{\sqrt{MT}}.
\]
Combining these estimates proves~\eqref{eq:localseg}.
\end{pf}

\paragraph{IC specialization.}
For the standard block Local SEG schedule, every communication round contains
$\tau:=\left\lfloor K/2\right\rfloor$
completed extragradient iterations followed by a terminal-state
synchronization. Hence $T=R\tau$. For $K\ge2$,
$\lfloor K/2\rfloor\ge K/3$, and Corollary~\ref{cor:localseg} gives
\begin{equation}
\label{eq:localseg-ic}
\Expct\!\left[\Gap_\Omega(\widetilde z_T)\right]
\le
C\left(
\frac{L\Omega^2}{KR}
+\frac{\sigma\Omega}{\sqrt{MKR}}
\right).
\end{equation}

\paragraph{Terminal-only aggregation.}
The bound~\eqref{eq:localseg-ic} is attained without intermediate
synchronizations. Divide the local computation into $R$ blocks of
$\lfloor K/2\rfloor$ SEG iterations, leave the first $R-1$ communication
opportunities unused, and aggregate the accumulated half-iterate sums at the
final opportunity. The resulting output is
\[
\widetilde z_T
=
\frac1M\sum_{m=1}^M
\left(
\frac1T\sum_{t=0}^{T-1}y_t^m
\right).
\]

Table~\ref{tab:main-summary} places the Local SEG upper bound beside two
standard SEG baselines, the Wang--Kale analysis of Local Extra SGD, and the two
lower bounds proved here.

\begin{table}[t]
\centering\footnotesize
\caption{Expected weak-gap upper and lower bounds in the IC model, up to
universal constants. The Local SEG upper bound and the two lower bounds in
bold are proved in this paper; in the lower-bound rows; $c>0$ denotes an
absolute numerical constant.}
\label{tab:main-summary}
\setlength{\tabcolsep}{2.5pt}
\renewcommand{\arraystretch}{1.10}
\begin{tabular}{@{}>{\raggedright\arraybackslash}p{0.23\linewidth}
                    >{\raggedright\arraybackslash}p{0.21\linewidth}
                    >{\raggedright\arraybackslash}p{0.34\linewidth}
                    >{\raggedright\arraybackslash}p{0.15\linewidth}@{}}
\toprule
\textbf{Method / result} & \textbf{Scope} & \textbf{Weak-gap bound}
& \textbf{Interpretation}\\
\midrule
\multicolumn{4}{@{}l}{\textbf{Upper bounds}}\\
\midrule
Minibatch SEG~\cite{Nemirovski2004,JNT2011}
& general monotone; bounded variance; $R\ge2$
& $\displaystyle \frac{L\Omega^2}{R}
  +\frac{\sigma\Omega}{\sqrt{MKR}}$
& deterministic horizon is $R$\\
\midrule
Single-device SEG~\cite{Korpelevich1976,Nemirovski2004,JNT2011}
& general monotone; bounded variance; $KR\ge2$
& $\displaystyle \frac{L\Omega^2}{KR}
  +\frac{\sigma\Omega}{\sqrt{KR}}$
& variance has single-device scaling\\
\midrule
Wang--Kale analysis of Local Extra SGD~\cite{WangKale2026}
& affine mean; conditionally unbiased bounded-variance oracle; $K\ge2$
& $\displaystyle \frac{L\Omega^2}{KR}
  +\frac{\sigma\Omega}{\sqrt{MKR}}
  +\frac{\Omega^{4/3}\sigma^{2/3}L^{1/3}}
  {K^{1/3}R^{2/3}}$
& proved bound contains a client-drift term\\
\midrule
\textbf{Local SEG upper bound (this paper)}
(Corollary~\ref{cor:localseg})
& affine monotone mean; conditional independence across devices within each
stage; $K\ge2$
& $\displaystyle \frac{L\Omega^2}{KR}
  +\frac{\sigma\Omega}{\sqrt{MKR}}$
& two-term rate matched in the stated class\\
\midrule
\multicolumn{4}{@{}l}{\textbf{Lower bounds proved in this paper}}\\
\midrule
\textbf{Affine Gaussian lower bound (this paper)}
(Theorem~\ref{thm:affinelb})
& unconstrained affine VI; Gaussian noise; zero-respecting depth
& $\displaystyle \ge c\left[
  \frac{L\Omega^2}{KR}
  +\min\!\left\{L\Omega^2,
  \frac{\sigma\Omega}{\sqrt{MKR}}\right\}
  \right]$
& order-matching before statistical saturation\\
\midrule
\textbf{Communication-locality lower bound (this paper)}
(Theorem~\ref{thm:communication-boundary})
& bounded affine VI; query-dependent noise; $\sigma>0$; feasible zero-respecting methods
& $\displaystyle \ge c\,\min\!\left\{
  \frac{L\Omega^2}{R\log(\mathrm e M)},
  \frac{\sigma^2}{LKR}
  \right\}$
& communication-dependent constrained regime\\
\bottomrule
\end{tabular}
\end{table}

The first four data rows are upper bounds, and the last two are lower bounds.
The upper-bound comparison uses homogeneous mean operators and unconstrained
Euclidean geometry. For our Local SEG bound, oracle errors are conditionally
centered and conditionally independent across devices within each stage. The
Wang--Kale and present Local SEG rows assume full participation, equal device
weights, and a common deterministic stepsize. The final row describes a
bounded-domain boundary case.

In the first five data rows, $\Omega$ is the radius of the root-centered
comparison ball; in the last row, it bounds the diameter of the feasible set.

\subsection{Matching affine lower bound}

\begin{fed}[Affine Gaussian instance]
\label{def:affine-gaussian}
For $L,\Omega>0$ and $\sigma\ge0$, an affine Gaussian instance on $\R^n$
consists of a mean operator $G$ satisfying
Assumptions~\ref{as:monotone}, \ref{as:lipschitz}, and~\ref{as:affine}, with
Lipschitz constant at most $L$, and a root $z^\star$ satisfying
$\norm{z^\star}\le\Omega.$
At every oracle call, each device receives
\[
\widehat G(z)=G(z)+\xi,
\]
where the errors are independent draws from the same centered, possibly
degenerate Gaussian distribution and satisfy
\[
\Expct\!\left[\norm{\xi}^2\right]\le\sigma^2.
\]
The notation $\Gap_\Omega$ refers to the root-centered comparison ball
in~\eqref{eq:gap-omega}. All expectations include both oracle and algorithmic
randomness.
\end{fed}

\begin{teo}[Affine Gaussian lower bound]
\label{thm:affinelb}
Let $L,\Omega>0$, $\sigma\ge0$, and let $M,K,R\ge1$ be integers. For every
$n\ge4KR+4$ and every distributed zero-respecting IC algorithm on $\R^n$,
there is an affine Gaussian instance in the sense of
Definition~\ref{def:affine-gaussian} such that
\begin{equation}
\label{eq:affine-lower}
\Expct\!\left[\Gap_\Omega(\hz)\right]\ge c\left[
\frac{L\Omega^2}{KR}
+\min\left\{L\Omega^2,
\frac{\sigma\Omega}{\sqrt{MKR}}\right\}\right].
\end{equation}
where $c>0$ is an absolute numerical constant independent of the problem
parameters and the algorithm.

For every $n\ge1$ and every randomized IC algorithm on
$\R^n$ with the same budgets, there is a possibly different affine Gaussian
instance with these parameters for which
\[
\Expct\!\left[\Gap_\Omega(\hz)\right]
\ge c\min\left\{L\Omega^2,
\frac{\sigma\Omega}{\sqrt{MKR}}\right\}.
\]
\end{teo}

\noindent\textit{Proof sketch.}
The deterministic and statistical terms are obtained from two affine
instances. For the deterministic term, a skew-symmetric bilinear zero chain
reveals at most one new coordinate per sequential oracle call. After depth
$KR$, a terminal coordinate remains hidden, and direct evaluation of the weak
gap gives a lower bound of order $L\Omega^2/(KR)$. For the statistical term, a
one-dimensional affine root observed through additive Gaussian noise yields a
two-point testing problem with $MKR$ observations. This gives the saturated
statistical term for every randomized IC algorithm. The maximum of the two
risks is at least half their sum, which proves~\eqref{eq:affine-lower}.
Appendix~\ref{app:lower} contains the complete argument.

The deterministic construction starts at the origin and has
$\norm{z^\star}\le\Omega$, so it belongs to the upper class with common
initialization $z_0=0$. Therefore the $KR$ and $MKR$ dependences in
Corollary~\ref{cor:localseg} are order-optimal in the stated zero-respecting
class whenever the statistical term is nonsaturated.

For the comparison in Table~\ref{tab:main-summary}, Wang and Kale count
extragradient iterations, whereas $K$ here counts oracle calls. Their
comparison set is a ball centered at the initial device mean. Taking its
radius to be $2\Omega$ covers the root-centered ball in
\eqref{eq:gap-omega}; these conversions affect only absolute constants.

\section{Limits of exact averaging}
\label{sec:limits}

This section identifies a boundary of the unconstrained exact-reduction
result. We retain an affine mean operator but require feasible queries and
outputs on a bounded domain.

\subsection{A communication-locality boundary}

The exact reduction and its two-resource rate concern unconstrained Euclidean
updates.
The next theorem considers a bounded feasible set, where the algorithm must
keep its queries and output feasible. For this result only, write
\begin{equation}
\label{eq:gap-feasible}
\Gap_{\Zset}(x):=\sup_{z\in\Zset}\inner{G(z)}{x-z},
\qquad x\in\Zset.
\end{equation}
This gap is nonnegative because the supremum includes $z=x$.

\paragraph{Scope.}
The exact identity applies to an affine common mean operator, unprojected
Euclidean updates, equal aggregation weights, a common deterministic stepsize,
and homogeneous device operators. Individual oracle realizations may be
nonlinear and query dependent, and their effective-noise laws may depend on the
synchronization schedule. The boundary result below concerns feasible
distributed zero-respecting algorithms on a bounded domain.

For $L,\Omega,\sigma>0$ and an integer $M\ge1$, set
\[
\ell_M:=\log(\mathrm e M),
\qquad
p:=\min\left\{1,\frac{L^2\Omega^2}{32\sigma^2}\right\}.
\]

\begin{teo}[Communication-locality boundary]
\label{thm:communication-boundary}
Let $L,\Omega,\sigma>0$, let $M,K,R\ge1$ be integers, and let $n$ be an
integer satisfying
$n\ge
2\left\lceil4R\bigl(Kp+\ell_M\bigr)\right\rceil+2.$
For every distributed zero-respecting IC algorithm on $\R^n$ that makes
feasible queries and returns a feasible output, there exist
a nonempty compact convex set
$\Zset\subseteq\R^n$ of diameter at most $\Omega$, an affine monotone
$L$-Lipschitz operator $G:\Zset\to\R^n$ with a root in $\Zset$, and a
query-dependent stochastic oracle satisfying Assumption~\ref{as:oracle} with
conditional variance bounded by $\sigma^2$. The oracle law is common to all
devices, and the output $\hz$ satisfies
\begin{equation}
\label{eq:communication-boundary}
\Expct\!\left[\Gap_{\Zset}(\hz)\right]
\ge c\min\left\{
\frac{L\Omega^2}{R\ell_M},
\frac{\sigma^2}{LKR}
\right\}.
\end{equation}
where $c>0$ is an absolute numerical constant.
\end{teo}

\noindent\textit{Proof sketch.}
A query-dependent Bernoulli oracle reveals the next zero-chain coordinate only
with probability $p$. Rescaling a rare reveal by $1/p$ preserves conditional
unbiasedness, while the bounded feasible geometry controls its variance. In
one IC round, the furthest of $M$ devices advances by only
$O(Kp+\log(\mathrm eM))$ coordinates in expectation. Choosing a chain longer
than the progress attainable in $R$ rounds leaves a hidden terminal coordinate
and hence a positive weak gap. Balancing $p$ with the variance budget gives
the two regimes in~\eqref{eq:communication-boundary}. Appendix
\ref{app:communication-boundary} supplies the full argument, which adapts the
robust zero-chain construction of Woodworth et al.~\cite{WBSS2021} to an
affine bilinear VI.

In the regime $\sigma^2\ge L^2\Omega^2K/\ell_M$, the first term in
\eqref{eq:communication-boundary} is active, so the error has an explicit
$1/R$ communication dependence. The two results cover distinct geometries.
Theorem~\ref{thm:exact} concerns
unprojected updates on $\R^n$, while the hard instance above combines active
feasibility geometry with a query-dependent revealing oracle.

\section{Numerical illustration of the two-resource effect}
\label{sec:numerics}

The experiment isolates the two resources in
Corollary~\ref{cor:localseg}: adaptive sequential depth and aggregate sampling
across devices. Up to integer floors, minibatch SEG, single-device SEG, and
Local SEG have the respective resource pairs
$(R,MKR)$, $(KR,KR)$, and $(KR,MKR)$, where the first entry is adaptive
depth and the second is the sample budget entering variance reduction. Thus
the two baselines separate the two advantages combined by Local SEG.

We use the ten-dimensional affine VI
\begin{equation}
\label{eq:numerical-affine}
G(z)=Jz,
\qquad
J=\frac{1}{\left\|\left[\begin{smallmatrix}0&A^\top\\-A&0\end{smallmatrix}\right]\right\|_{\rm op}}
\begin{bmatrix}0&A^\top\\-A&0\end{bmatrix},
\end{equation}
where $A\in\R^{5\times5}$ is one fixed Gaussian draw. Thus $J^\top=-J$,
$\norm{J}_{\rm op}=1$, and $z^\star=0$. The comparison set is the unit ball,
and a call returns $Jz+\xi$ with
$\xi\sim\mathcal N(0,I_{10}/10)$ independently across devices and stages.
For a comparison radius $\Omega$, the restricted weak gap is
$\Gap_\Omega(x)=\Omega\norm{J^\top x}$; the experiment uses $\Omega=1$.

\begin{figure}[h]
\centering
\includegraphics[width=0.84\textwidth]{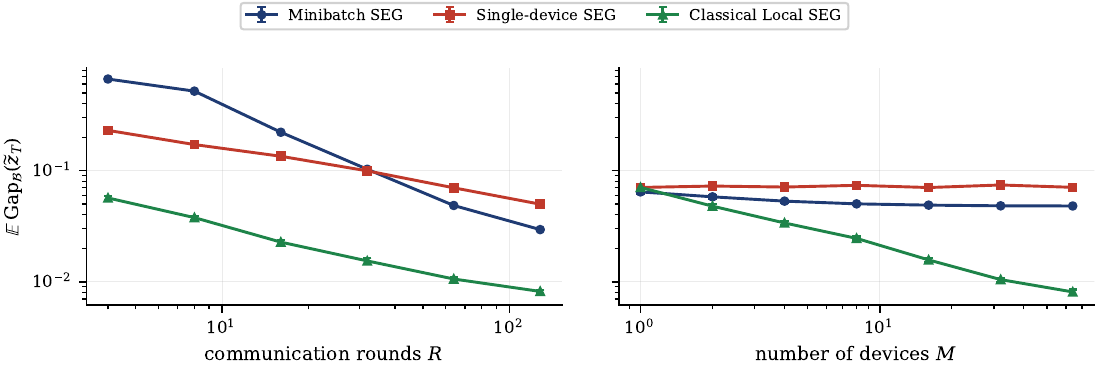}
\caption{Separation of sequential-depth and device-averaging effects on the
affine problem~\eqref{eq:numerical-affine}. Left: expected restricted weak gap
as the communication budget increases. Right: expected restricted weak gap as
the number of devices increases. Error bars are pointwise normal-approximation
95\% Monte Carlo confidence intervals from 80 independent replications.}
\label{fig:resource-scaling}
\end{figure}

Figure~\ref{fig:resource-scaling} compares minibatch SEG, single-device SEG,
and classical Local SEG. Every method uses the same analytical,
data-independent stepsize rule and the exact two-calls-per-SEG accounting described in
Appendix~\ref{app:numerics}. The left panel fixes $M=32$ and $K=16$ and varies
the communication budget. The right panel fixes $K=16$ and $R=64$ and varies
the number of devices.

The left panel isolates adaptive depth: at the same communication and
aggregate-oracle budgets, Local SEG extends the horizon from order $R$ to
order $KR$ and improves on minibatch SEG.

The right panel isolates device averaging. The single-device curve is nearly
insensitive to $M$. Minibatch SEG benefits from larger batches but approaches
the deterministic floor associated with its shorter horizon, whereas Local SEG
retains the longer horizon while its effective variance decreases with $M$.
This behavior is consistent with the deterministic and stochastic terms in
\eqref{eq:localseg-ic}.

\section{Conclusion}
\label{sec:conclusion}

We studied the classical Local SEG method for stochastic affine monotone VIs
under intermittent communication. For unconstrained Euclidean updates with
equal device weights and a common stepsize, affinity makes device averaging
commute with both stages of the extragradient recursion. The identity is
pathwise and continues to hold when individual oracle realizations are
nonlinear and query dependent. Under conditionally centered errors that are
conditionally independent across devices within each stage, the averaged
oracle error has conditional second moment at most $\sigma^2/M$. For $K\ge2$,
this gives
\[
\mathcal O\!\left(
\frac{L\Omega^2}{KR}
+\frac{\sigma\Omega}{\sqrt{MKR}}
\right)
\]
for the expected restricted weak gap. The deterministic component depends on
the local sequential depth, whereas the stochastic component depends on the
aggregate oracle budget. The same finite-horizon guarantee is attained with
one terminal aggregation.

The affine Gaussian lower bound matches these resource dependences before
statistical saturation for distributed zero-respecting algorithms, while its
statistical component applies to arbitrary randomized IC algorithms. The
bounded-domain construction gives a complementary communication-locality lower
bound for feasible distributed zero-respecting methods. Together, the results
separate the unconstrained affine regime, in which the device-mean recursion
closes exactly, from a constrained regime in which active geometry and
query-dependent noise can restrict the propagation of information.

Extensions to approximately affine operators, projected or mirror updates,
unequal aggregation weights, and heterogeneous device operators remain open.
Interpolation-type noise models constitute a further direction because the
effective noise may then depend on disagreement among local query points.

\section*{Acknowledgments}
This work was supported by Ministry of Science and Higher Education of
the Russian Federation (Grant No. 075-15-2024-529). 

\noindent We thank Georgy Bondarenko for comments and discussions. 

\paragraph{Declaration on the use of generative AI.}
Generative AI tools were used in a supporting role for language editing of the final manuscript, presentation-level organization, and limited assistance with research code. The research questions, mathematical results, experimental design, and scientific conclusions were developed by the authors. All AI-assisted text and code, as well as all proofs and numerical results, were reviewed and independently verified by the authors, who take full responsibility for the content of the manuscript.

\newpage
\appendix

\section{Affine lower-bound constructions}
\label{app:lower}

This appendix proves the two lower-bound theorems. A bilinear zero chain gives
the deterministic depth obstruction, and a one-dimensional two-point
experiment gives the Gaussian statistical obstruction. The same chain,
combined with a rare-revelation oracle, proves the communication-locality
boundary on a bounded domain.

For an operator $H:\mathcal D\to\R^d$, a nonempty comparison set
$\mathcal C\subseteq\mathcal D$, and $x\in\R^d$, define
\[
\Gap_{\mathcal C}^{H}(x):=
\mathop{\rm sup}_{z\in\mathcal C}\inner{H(z)}{x-z}.
\]
When the operator is clear, we omit the superscript. Thus
$\Gap_{\mathcal B}$ and $\Gap_{\mathcal U\times\mathcal V}$ in this appendix
are concrete instances of the same restricted weak-gap definition as
$\Gap_\Omega$ in~\eqref{eq:gap-omega} and $\Gap_{\Zset}$
in~\eqref{eq:gap-feasible}.

\subsection{Formal IC transcript and support condition}
\label{app:ic-formal}

For the lower bound it is enough to consider the following maximal-transcript
envelope, which permits full sharing after every completed round. Let
$x^m_{k,r}$ denote the $k$th query of device $m$ in round $r$, and let
$U_{\rm alg}$ collect the shared algorithmic randomness.
The query $x^m_{k,r}$ is measurable with respect to $U_{\rm alg}$, all
query--reply pairs from rounds $r'<r$, and the same device's query--reply pairs
with indices $k'<k$ in round $r$. The output $\hz$ is measurable with respect
to the full transcript and $U_{\rm alg}$. A distributed zero-respecting
algorithm additionally satisfies, almost surely,
\begin{align}
\label{eq:zero-respecting-query}
\supp(x^m_{k,r})
&\subseteq
\bigcup_{\substack{1\le m'\le M,\ 1\le k'\le K\\1\le r'<r}}
\supp\!\left(\widehat G(x^{m'}_{k',r'},\zeta^{m'}_{k',r'})\right)
\cup
\bigcup_{1\le k'<k}
\supp\!\left(\widehat G(x^m_{k',r},\zeta^m_{k',r})\right),\\
\label{eq:zero-respecting-output}
\supp(\hz)
&\subseteq
\bigcup_{\substack{1\le m\le M,\ 1\le k\le K\\1\le r\le R}}
\supp\!\left(\widehat G(x^m_{k,r},\zeta^m_{k,r})\right).
\end{align}
Shared randomness may choose measurable coefficients subject to
\eqref{eq:zero-respecting-query}--\eqref{eq:zero-respecting-output}. The
statistical construction below applies to arbitrary randomized IC algorithms.

\subsection{Bilinear zero chain}

\begin{lem}[Bilinear chain]
\label{lem:chain}
Let $N\ge1$ be an integer and let $r,\kappa>0$. Define the upper bidiagonal
matrix $A_N\in\R^{N\times N}$ by
$(A_Nu)_i=u_i-u_{i+1}$ for $i<N$ and $(A_Nu)_N=u_N$. Put
$\eta=r/\sqrt N$ and
\[
\Phi(u,v):=\kappa\bigl(\inner{A_Nu}{v}+\eta u_1\bigr),
\qquad
G_N(u,v):=\kappa\bigl(A_N^\top v+\eta e_1,-A_Nu\bigr).
\]
Then $G_N$ is affine, monotone, and $2\kappa$-Lipschitz on $\R^{2N}$,
with root $(u^\star,v^\star)=(0,-\eta\mathbf 1)$. Let
$\mathcal U=\mathcal V:=\{w\in\R^N:\norm w\le r\}$. For every
$(\hat u,\hat v)\in\R^{2N}$ with $\hat v_N=0$,
\begin{equation}
\label{eq:chain-gap}
\Gap_{\mathcal U\times\mathcal V}(\hat u,\hat v)
\ge \frac{\kappa r^2}{N}.
\end{equation}
Moreover, in the interleaved ordering
$(u_1,v_1,u_2,v_2,\ldots)$, if $w\in\R^{2N}$ is supported on the first $j$
coordinates, where $0\le j\le2N$, then $G_N(w)$ is supported on the first
$\min\{j+1,2N\}$ coordinates.
\end{lem}

\begin{pf}
The linear part of $G_N$ is skew-symmetric, so the operator is monotone. Since
$\norm{A_N}_{\rm op}\le2$, its Lipschitz constant is at most $2\kappa$.
Solving $A_Nu=0$ and $A_N^\top v=-\eta e_1$ gives
$u^\star=0$ and $v^\star=-\eta\mathbf 1$.

For any candidate $(\hat u,\hat v)$, direct expansion of the weak gap gives
\begin{align*}
\Gap_{\mathcal U\times\mathcal V}(\hat u,\hat v)
&=\sup_{u\in\mathcal U,\,v\in\mathcal V}
\bigl[\Phi(\hat u,v)-\Phi(u,\hat v)\bigr]=\sup_{v\in\mathcal V}\Phi(\hat u,v)
-\inf_{u\in\mathcal U}\Phi(u,\hat v).
\end{align*}
Because $\sum_i(A_N\hat u)_i=\hat u_1$,
\[
\sup_{v\in\mathcal V}\Phi(\hat u,v)
=\kappa\eta\hat u_1+\kappa r\norm{A_N\hat u}\ge0.
\]
If $\hat v_N=0$, the entries of
$A_N^\top\hat v+\eta e_1$ sum to $\eta$, and hence
\[
-\inf_{u\in\mathcal U}\Phi(u,\hat v)
=\kappa r\norm{A_N^\top\hat v+\eta e_1}
\ge\frac{\kappa r\eta}{\sqrt N}
=\frac{\kappa r^2}{N}.
\]
This proves~\eqref{eq:chain-gap}. For the support claim, set $v_0=0$ and
$u_{N+1}=0$. In the interleaved ordering,
\[
(G_N)_{u_i}
=\kappa\bigl(v_i-v_{i-1}+\eta\mathbf 1_{\{i=1\}}\bigr),
\qquad
(G_N)_{v_i}=\kappa(u_{i+1}-u_i).
\]
These coordinate identities prove the stated support property.
\end{pf}

\subsection{Statistical affine-root obstruction}

\begin{lem}[Affine-root statistical lower bound]
\label{lem:statlbroot}
Let $L,\Omega>0$, $\sigma\ge0$, and let $N_{\rm q}\ge1$ be an integer. For
every possibly randomized adaptive algorithm using at most $N_{\rm q}$ oracle
calls, there is a one-dimensional affine monotone operator with Lipschitz
constant at most $L$, a root $z^\star$ satisfying
$|z^\star|\le\Omega/4$, and an additive Gaussian oracle with independent
noise of variance at most $\sigma^2$ such that, for
$\mathcal B:=\{z:|z-z^\star|\le\Omega\}$,
\begin{equation}
\label{eq:statistical-root-lower}
\Expct\!\left[\Gap_{\mathcal B}(\hz)\right]
\ge c\min\left\{L\Omega^2,
\frac{\sigma\Omega}{\sqrt{N_{\rm q}}}\right\}.
\end{equation}
where $c>0$ is an absolute numerical constant.
\end{lem}

\begin{pf}
The claim is immediate when $\sigma=0$, so suppose $\sigma>0$. Let
$\theta\in\{-\beta,+\beta\}$ with $\beta=\Omega/4$ and set
\[
G_\theta(z):=a(z-\theta),
\qquad
a:=\min\left\{L,\frac{c_0\sigma}{\Omega\sqrt{N_{\rm q}}}\right\},
\]
where $c_0>0$ is a sufficiently small numerical constant. The oracle returns
$G_\theta(z)+\xi$ with $\xi\sim\mathcal N(0,\sigma^2)$ independently at each
call. At every query, the two conditional oracle means differ by $2a\beta$.
The chain rule for relative entropy therefore bounds the KL divergence between
the two adaptive transcripts by
$N_{\rm q}(2a\beta)^2/(2\sigma^2)$. Le Cam's two-point lemma
\cite[Section~2.2]{Tsybakov2009} then gives a constant lower bound on the
probability of confusing the two roots.

For $x\in\R$, write $e=x-\theta$ and $q=z-\theta$. Directly from the weak-gap
definition,
\[
\Gap_{\mathcal B}(x)
=a\sup_{|q|\le\Omega}q(e-q)
\ge c\,a\min\{e^2,\Omega^2\}.
\]
Indeed, take $q=e/2$ when $|e|\le2\Omega$ and
$q=\Omega\operatorname{sign}(e)$ otherwise. Predicting the sign of $\theta$
from the sign of $\hz$, a testing error implies
$|\hz-\theta|\ge\Omega/4$. Thus the gap on that event is at least
$ca\Omega^2$, and the testing lower bound yields
\eqref{eq:statistical-root-lower}.
More explicitly, averaging over the two equiprobable signs of $\theta$ gives
this lower bound for the average risk, so at least one of the two operators
satisfies~\eqref{eq:statistical-root-lower}.
\end{pf}

\begin{pf}[Proof of Theorem~\ref{thm:affinelb}]
For the deterministic term, take the chain in Lemma~\ref{lem:chain} with
\[
N_0:=2KR+2,
\qquad
\kappa:=L/2,
\qquad
r:=\Omega/\sqrt5,
\]
and use the noiseless, hence degenerate Gaussian, oracle. At the start of round
$s$, suppose the shared transcript is supported on the first $(s-1)K$
interleaved coordinates. By~\eqref{eq:zero-respecting-query}, the first query
of every device is supported on this prefix. Lemma~\ref{lem:chain} permits at
most one new coordinate per local reply, so after its $k$th query a device has
revealed at most the first $(s-1)K+k$ coordinates. The union over devices
remains supported on the same prefix. Induction over $k$ and $s$ gives final
progress at most $KR$; in particular, the output satisfies $\hat v_{N_0}=0$.
The product set
$\mathcal U\times\mathcal V$ lies in the ball of radius $\Omega$ centered at
$z^\star=(0,-\eta\mathbf1)$, because
\[
\norm{(u,v)-z^\star}^2
\le r^2+(\norm v+\norm{v^\star})^2
\le5r^2=\Omega^2.
\]
Consequently, Lemma~\ref{lem:chain} gives
\[
\Gap_\Omega(\hz)
\ge\frac{\kappa r^2}{N_0}
\ge\frac{L\Omega^2}{40KR}.
\]
The active dimension is $2N_0=4KR+4$; extra coordinates are padded with
zeros.

For the statistical term, apply Lemma~\ref{lem:statlbroot} with
$N_{\rm q}=MKR$, place the scalar construction in the first coordinate, and
set all remaining operator and noise coordinates to zero. Its set $\mathcal B$
embeds in the first coordinate of the root-centered ball defining
$\Gap_\Omega$. Because the padded operator vanishes in all other coordinates,
the resulting gap equals the scalar gap. Selecting the hard instance with the
larger of the deterministic and statistical risks and
using $\max\{a,b\}\ge(a+b)/2$ proves~\eqref{eq:affine-lower}. The statistical
construction applies to arbitrary randomized IC algorithms, which proves the
second assertion of the theorem.
\end{pf}

\subsection{Communication-locality construction}
\label{app:communication-boundary}

\begin{pf}[Proof of Theorem~\ref{thm:communication-boundary}]
We use Lemma~\ref{lem:chain} with
\[
\kappa:=L/2,
\qquad
r:=\frac{\Omega}{2\sqrt2}.
\]
For any chain length $N$, the product domain
$\Zset:=\mathcal U\times\mathcal V$ has diameter $2\sqrt2r=\Omega$.
Its root $(0,-\eta\mathbf1)$ is feasible because
$\norm{\eta\mathbf1}=r$, and $G_N$ is $L$-Lipschitz.

We first define the stochastic oracle. Write the coordinates in the
interleaved order from Lemma~\ref{lem:chain}, and put $d:=2N$. For
$0\le j\le d$, let $P_j$ denote the orthogonal projection onto the first $j$
interleaved coordinates, with $P_0=0$ and $P_d=I$. For $w\in\Zset$, define
\[
\operatorname{prog}(w):=\max\{j:w_j\ne0\},
\qquad \operatorname{prog}(0):=0.
\]
With $j=\operatorname{prog}(w)$, set
\[
H_0(w):=P_jG_N(w),
\qquad
H_1(w):=G_N(w),
\qquad
D(w):=H_1(w)-H_0(w).
\]
The zero-chain formulas in Lemma~\ref{lem:chain} show that $D(w)$ contains
only the possible next frontier coordinate. More explicitly,
letting $e_{u_i}$ and $e_{v_i}$ denote the basis vectors of the corresponding
interleaved coordinates,
\[
D(w)=
\begin{cases}
\kappa\eta e_{u_1}, & j=0,\\
-\kappa u_i e_{v_i}, & j=2i-1,\\
-\kappa v_i e_{u_{i+1}}, & j=2i<2N,\\
0, & j=2N.
\end{cases}
\]
Since $w=(u,v)\in\mathcal U\times\mathcal V$ and
$\eta\le r$, this gives
\begin{equation}
\label{eq:frontier-bound}
\norm{D(w)}\le\kappa r,
\qquad
\kappa^2r^2=\frac{L^2\Omega^2}{32}.
\end{equation}

At each call, draw a new $U\sim\operatorname{Bernoulli}(p)$, independently
across devices and calls conditional on the past, and return
\begin{equation}
\label{eq:rare-reveal-oracle}
\widehat G_N(w;U):=H_0(w)+\frac{U}{p}D(w).
\end{equation}
Because $\sigma>0$, the probability $p$ belongs to $(0,1]$. Conditional on the
pre-call history, $w$ is fixed and
\begin{align}
\label{eq:rare-reveal-moments}
\Expct\!\left[\widehat G_N(w;U)\mid\mathcal F\right]
&=G_N(w),\\
\Expct\!\left[
\norm{\widehat G_N(w;U)-G_N(w)}^2\mid\mathcal F
\right]
&=\frac{1-p}{p}\norm{D(w)}^2\le\sigma^2.
\end{align}
Indeed, the last inequality is immediate for $p=1$. If $p<1$, then
$p=L^2\Omega^2/(32\sigma^2)$, and~\eqref{eq:frontier-bound} gives the claim.
Thus the same affine mean operator and the same oracle kernel are used on all
devices, although the oracle realization depends on the query through its
revealed support.

It remains to bound the revealed progress. Let $\mathcal R_s$ be the finite
set of oracle replies produced through round $s$, and define
\[
S_s:=\max\left(\{0\}\cup
\bigcup_{g\in\mathcal R_s}\supp(g)\right).
\]
Before the first reply, $S_0=0$. A call with $U=1$ reveals the next coordinate.
Let $\mathcal H_{s-1}$ denote the history before round $s$. Conditional on
$\mathcal H_{s-1}$, complete the $K$-call sequence with independent coins and
set
\[
X_{m,s}:=\sum_{k=1}^K U_{m,k,s}
\sim\operatorname{Binomial}(K,p).
\]
The support condition~\eqref{eq:zero-respecting-query} then gives, almost
surely,
\begin{equation}
\label{eq:round-progress}
S_s-S_{s-1}\le\max_{1\le m\le M}X_{m,s}.
\end{equation}
For completeness, the conditional moment-generating-function estimate
\[
\Expct\!\left[
\exp\!\left(\max_m X_{m,s}\right)\mid\mathcal H_{s-1}
\right]
\le M(1-p+p\mathrm e)^K
\le\exp\!\left(\log M+(\mathrm e-1)Kp\right)
\]
and conditional Jensen's inequality imply
\[
\Expct\!\left[\max_m X_{m,s}\mid\mathcal H_{s-1}\right]
\le2\bigl(Kp+\ell_M\bigr).
\]
Iterating~\eqref{eq:round-progress} yields
\begin{equation}
\label{eq:expected-progress}
\Expct\!\left[S_R\right]
\le2R\bigl(Kp+\ell_M\bigr).
\end{equation}

Choose the active chain length
\[
N:=\left\lceil4R\bigl(Kp+\ell_M\bigr)\right\rceil+1.
\]
Markov's inequality and~\eqref{eq:expected-progress} give
$\mathbb P\{S_R<N\}\ge1/2$. On this event,
\eqref{eq:zero-respecting-output} implies that the output has
$\hat v_N=0$. Lemma~\ref{lem:chain} therefore gives
$\Gap_{\Zset}(\hz)\ge\kappa r^2/N$. Off this event the gap remains
nonnegative because $\hz$ is feasible. Consequently,
\begin{equation}
\label{eq:raw-communication-lower}
\Expct\!\left[\Gap_{\Zset}(\hz)\right]
\ge c_0\frac{L\Omega^2}{R(Kp+\ell_M)}
\end{equation}
where $c_0>0$ is an absolute numerical constant.

To express the two regimes, put $q:=\sigma^2/(L^2\Omega^2)$. If $p=1$, then
$q\le1/32$ and
\[
\frac1{K+\ell_M}
\ge\frac12\min\left\{\frac1{\ell_M},\frac qK\right\}.
\]
If $p<1$, then $p=1/(32q)$, and the same inequality, with a possibly smaller
absolute constant, follows from
\[
\frac1{Kp+\ell_M}
=\frac1{K/(32q)+\ell_M}
\ge\frac12\min\left\{\frac{32q}{K},\frac1{\ell_M}\right\}
\ge\frac12\min\left\{\frac qK,\frac1{\ell_M}\right\}.
\]
Multiplying by $L\Omega^2/R$ in
\eqref{eq:raw-communication-lower} proves
\eqref{eq:communication-boundary}.

Finally, the dimension condition in the theorem is exactly $n\ge2N$.
Embed the product balls in the first $2N$ coordinates and pad the operator
and oracle replies by zeros when $n>2N$.
\end{pf}

\section{Unprojected stochastic extragradient estimate}
\label{app:seg}

The main upper bound follows from a single centralized SEG estimate. Its
adapted-error formulation accommodates the averaged residuals in
Theorem~\ref{thm:exact}, whose distribution may differ from that of a primitive
oracle call at the mean query.

\begin{pro}[Unprojected SEG with adapted errors]
\label{prop:freeeg}
Let $G:\R^n\to\R^n$ be monotone and $L$-Lipschitz for some $L>0$, and
suppose $G(z^\star)=0$. Fix an integer $T\ge1$, a stepsize
$0<\gamma\le1/(2L)$, and $v\ge0$. Let
$\mathcal F_0\subseteq\mathcal F_{1/2}\subseteq\mathcal F_1\subseteq\cdots
\subseteq\mathcal F_{T-1/2}\subseteq\mathcal F_T$ be a filtration. For
$t=0,\ldots,T-1$, let $z_t$ be
$\mathcal F_t$-measurable. Let
$\Delta_t$ and $\Delta_{t+1/2}$ be measurable with respect to
$\mathcal F_{t+1/2}$ and $\mathcal F_{t+1}$, respectively, and suppose
\begin{align}
\label{eq:adapted-errors}
\Expct\!\left[\norm{\Delta_t}^2\mid\mathcal F_t\right]&\le v,\\
\Expct\!\left[\Delta_{t+1/2}\mid\mathcal F_{t+1/2}\right]&=0,
&\Expct\!\left[\norm{\Delta_{t+1/2}}^2
\mid\mathcal F_{t+1/2}\right]&\le v
\end{align}
Run
\[
y_t=z_t-\gamma\bigl(G(z_t)+\Delta_t\bigr),
\qquad
z_{t+1}=z_t-\gamma\bigl(G(y_t)+\Delta_{t+1/2}\bigr).
\]
If, for some $\Omega>0$,
$\norm{z_0-z^\star}\le\Omega$ almost surely and
$\mathcal B:=\{u:\norm{u-z^\star}\le\Omega\}$, then, for
$\widetilde z_T:=T^{-1}\sum_{t=0}^{T-1}y_t$,
\begin{equation}
\label{eq:freeeg}
\Expct\!\left[\Gap_{\mathcal B}(\widetilde z_T)\right]
\le\frac{3\Omega^2}{\gamma T}+9\gamma v.
\end{equation}
\end{pro}

\begin{pf}
Fix $u\in\mathcal B$. Expanding the two unprojected updates gives
\begin{align*}
2\gamma\inner{G(y_t)}{y_t-u}
&=\norm{z_t-u}^2-\norm{z_{t+1}-u}^2
-\norm{y_t-z_t}^2-\norm{y_t-z_{t+1}}^2\\
&\quad+2\gamma\inner{G(y_t)-G(z_t)}{y_t-z_{t+1}}
-2\gamma\inner{\Delta_t}{y_t-z_{t+1}}\\
&\quad+2\gamma\inner{\Delta_{t+1/2}}{z_t-z_{t+1}}
+2\gamma\inner{\Delta_{t+1/2}}{u-z_t}.
\end{align*}
Put $a=\norm{y_t-z_t}$, $b=\norm{y_t-z_{t+1}}$,
$d=\norm{\Delta_t}$, and $e=\norm{\Delta_{t+1/2}}$. Since
$\gamma L\le1/2$, Lipschitz continuity and Young's inequality give
\[
-a^2-b^2+2\gamma Lab+2\gamma db+2\gamma e(a+b)
\le8\gamma^2(d^2+e^2).
\]
Indeed, the four positive terms are bounded by
$(a^2+b^2)/2$, $b^2/4+4\gamma^2d^2$,
$a^2/4+4\gamma^2e^2$, and $b^2/4+4\gamma^2e^2$, respectively.
Substitution into the preceding identity leaves
\begin{align}
\label{eq:free-one-step}
2\gamma\inner{G(y_t)}{y_t-u}
&\le\norm{z_t-u}^2-\norm{z_{t+1}-u}^2\nonumber\\
&\quad+8\gamma^2
\bigl(\norm{\Delta_t}^2+\norm{\Delta_{t+1/2}}^2\bigr)
+2\gamma\inner{\Delta_{t+1/2}}{u-z_t}.
\end{align}

Monotonicity gives
$\inner{G(y_t)}{y_t-u}\ge\inner{G(u)}{y_t-u}$. To control the supremum over
$u$ prior to taking expectations, let $p_0=z^\star$ and let
$p_{t+1}$ be the Euclidean projection of
$p_t+\gamma\Delta_{t+1/2}$ onto $\mathcal B$. The projection inequality yields
\begin{equation}
\label{eq:online-comparator}
\sup_{u\in\mathcal B}\sum_{t=0}^{T-1}
\inner{\Delta_{t+1/2}}{u-p_t}
\le\frac{\Omega^2}{2\gamma}
+\frac{\gamma}{2}\sum_{t=0}^{T-1}\norm{\Delta_{t+1/2}}^2.
\end{equation}
The sequence $p_t-z_t$ is $\mathcal F_{t+1/2}$-measurable. Conditional
centering therefore gives
\[
\Expct\!\left[
\sum_{t=0}^{T-1}\inner{\Delta_{t+1/2}}{p_t-z_t}
\right]=0.
\]
Sum~\eqref{eq:free-one-step}, take the supremum over $u\in\mathcal B$, and
use~\eqref{eq:online-comparator}. Since
$\sup_{u\in\mathcal B}\norm{z_0-u}^2\le4\Omega^2$, the expected right-hand
side is at most $5\Omega^2+17\gamma^2Tv$. Dividing by $2\gamma T$ and rounding
the constants upward proves~\eqref{eq:freeeg}.
\end{pf}

\section{Experimental details}
\label{app:numerics}

This appendix records the exact accounting used
in Section~\ref{sec:numerics}. 

For a positive horizon $T$ and effective variance $V$, every method uses
\[
\gamma=\min\left\{\frac1{2L},\frac{\Omega}{\sqrt{TV}}\right\}.
\]
The second entry is interpreted as $+\infty$ when $V=0$.
The following table gives the values used by the implementation; a completed
SEG iteration always consumes two oracle calls.

\begin{table}[ht]
\centering\small
\caption{Iteration horizons, effective variances, and actual oracle calls.}
\label{tab:numerical-accounting}
\begin{tabular}{@{}lccc@{}}
\toprule
\textbf{Method} & $T$ & $V$ & \textbf{Actual calls}\\
\midrule
Minibatch SEG ($R\ge2$) & $\lfloor R/2\rfloor$ & $\sigma^2/(MK)$
& $2TMK$\\
\midrule
Single-device SEG ($KR\ge2$) & $\lfloor KR/2\rfloor$ & $\sigma^2$
& $2T$\\
\midrule
Classical Local SEG ($K\ge2$) & $R\lfloor K/2\rfloor$ & $\sigma^2/M$
& $2TM$\\
\bottomrule
\end{tabular}
\end{table}

The corresponding horizons are positive under $R\ge2$ for minibatch SEG,
$KR\ge2$ for single-device SEG, and $K\ge2$ for classical Local SEG. All
configurations in Figure~\ref{fig:resource-scaling} satisfy these conditions.

The experiment uses $L=\Omega=\sigma=1$ and the fixed pseudorandom seed
$20260815$. The matrix $A$ in~\eqref{eq:numerical-affine} and the common
initial point $z_0$ are drawn once and reused in every configuration; the
initial point is rescaled so that $\norm{z_0}=0.8$. The rounds grid is
$\{4,8,16,32,64,128\}$ and the devices grid is
$\{1,2,4,8,16,32,64\}$. The left panel fixes $M=32$ and $K=16$; the right
panel fixes $K=16$ and $R=64$. Since these values of $K$ and $R$ are even,
minibatch SEG and Local SEG use exactly $MKR$ oracle calls in the corresponding
comparisons.

Each point in Figure~\ref{fig:resource-scaling} is computed from 80 independent
replications. The error bars equal $1.96$ times the sample standard error. The
output is the temporal average of the half-iterates, and its merit value is the
closed-form gap following~\eqref{eq:numerical-affine}. The stepsize is fixed in
advance according to the displayed rule.

\small
\bibliographystyle{unsrt}
\bibliography{references}

\ifNDofficial
}
\label{article_end}
\endpaper
\fi

\end{document}